\documentclass{amsart}
\usepackage{amsrefs}
\usepackage{graphics}
\newcommand{\dontprint}[1]{\relax}
\title{The work of Andrei Okounkov}
\author{Giovanni Felder, ETH Zurich}
\begin{document}
\maketitle Andrei Okounkov's initial area of research is group
representation theory, with particular emphasis on combinatorial and
asymptotic aspects. He used this subject as a starting point to
obtain spectacular results in many different areas of mathematics
and mathematical physics, from complex and real algebraic geometry
to statistical mechanics, dynamical systems, probability theory and
topological string theory. The research of Okounkov has its roots in
very basic notions such as partitions, which form a recurrent theme
in his work. A partition $\lambda$ of a natural number $n$ is a
non-increasing sequence of integers
$\lambda_1\geq\lambda_2\geq\cdots\geq 0$ adding up to $n$.
Partitions are a basic combinatorial notion at the heart of the
representation theory. Okounkov started his career in this field in
Moscow where he worked with G. Olshanski, through whom he came in
contact with A. Vershik and his school in St.~Petersburg, in
particular S. Kerov. The research programme of these mathematicians,
to which Okounkov made substantial contributions, has at its core
the idea that partitions and other notions of representation theory
should be considered as random objects with respect to natural
probability measures. This idea was further developed by Okounkov,
who showed that, together with insights from geometry and ideas of
high energy physics, it can be applied to the most diverse areas of
mathematics.

This is an account of some of the highlights mostly of his recent
research.

I am grateful to Enrico Arbarello for explanations and for providing
me with very useful notes on Okounkov's work in algebraic geometry
and its context.

\section{Gromov--Witten invariants}

The context of several results of Okounkov and collaborators is the
theory of Gromov--Witten (GW) invariants. This section is a short
account of this theory. GW invariants originate from classical
questions of enumerative geometry, such as: how many rational curves
of degree $d$ in the plane go through $3d-1$ points in general
position? A completely new point of view on this kind of problems
appeared at the end of the eighties, when string theorists, working
on the idea that space-time is the product of four-dimensional
Minkowski space with a Ricci-flat compact complex three-fold, came
up with a prediction for the number of rational curves of given
degree in the quintic $x_1^5+\cdots+x_5^5=0$ in $\mathbb C P^4$.
Roughly speaking, physics gives predictions for differential
equations obeyed by generating functions of numbers of curves.
Solving these equations in power series gives recursion relations
for the numbers. In particular a recursion relation of Kontsevich
gave a complete answer to the above question on rational curves in
the plane.

In general, Gromov--Witten theory deals with intersection numbers on
moduli spaces of maps from curves to complex manifolds. Let $V$ be a
nonsingular projective variety over the complex numbers. Following
Kontsevich, the compact moduli space $\overline M_{g,n}(V,\beta)$ (a
Deligne--Mumford stack) of {\em stable maps} of class $\beta\in
H_2(V)$ is the space of isomorphism classes of data
$(C,p_1,\dots,p_n,f)$ where $C$ is a complex projective connected
nodal curve of genus $g$ with $n$ marked smooth points $p_1$, \dots,
$p_n$ and $f\colon C\to V$ is a stable map such that $[f(C)]=\beta$.
Stable means that if $f$ maps an irreducible component to a point
then this component should have a finite automorphism group. For
each $j=1,\dots,n$ two natural sets of cohomology classes can be
defined on these moduli space: 1) pull-backs
$\mathrm{ev}_j^*\alpha\in H^*(\overline M_{g,n}(V,\beta))$ of
cohomology classes $\alpha\in H^*(V)$ on the target $V$ by the
evaluation map $\mathrm{ev}_j\colon (C,p_1,\dots,p_n,f)\mapsto
f(p_j)$; 2) the powers of the first Chern class $\psi_j=c_1(L_j)\in
H^2(\overline M_{g,n}(V,\beta))$ of the line bundle $L_j$ whose
fiber at $(C,p_1,\dots,p_n,f)$ is the cotangent space $T^*_{p_j}C$
to $C$ at $p_j$. The {\it Gromov--Witten invariants} of $V$ are the
intersection numbers
\[
\langle \tau_{k_1}(\alpha_1)\cdots\tau_{k_n}(\alpha_n)\rangle^V_{\beta,g}
=\int_{\overline M_{g,n}(V,\beta)}\prod \psi_j^{k_j}\mathrm{ev}_j^*\alpha_j.
\]
If all $k_i$ are zero and the $\alpha_i$ are Poincar\'e duals of
subvarieties, the Gromov--Witten invariants have the interpretation
of counting the number of curves intersecting these subvarieties. As
indicated by Kontsevich, to define the integral one needs to
construct a {\em virtual fundamental class}, a homology class of
degree equal to the ``expected dimension''
\begin{equation}\label{e-1}
\mathrm{vir}\,\mathrm{dim}\,\overline M_{g,n}(V,\beta)
=-\beta\cdot K_V+(g-1)(3-\mathrm{dim}\,V)+n,
\end{equation}
where $K_V$ is the canonical class of $V$. This class was
constructed in works of Behrend--Fantechi and Li--Tian.

The theory of Gromov--Witten invariants is already non-trivial and
deep in the case where $V$ is a point. In this case $\overline
M_{g,n}=\overline M_{g,n}(\{\mathrm{pt}\})$ is the Deligne--Mumford
moduli space of stable curves of genus $g$ with $n$ marked points.
Witten conjectured and Kontsevich proved that the generating
function
\[
F(t_0,t_1,\dots)=\sum_{n=0}^\infty \frac1{n!}
\sum_{k_1+\cdots+k_n=3g-3+n}\langle\tau_{k_1}\cdots \tau_{k_n}\rangle^{\mathrm{pt}}_{g}\prod_{j=1}^n
{t_{k_j}},
\]
involving simultaneously all genera and numbers of marked points,
obeys an infinite set of partial differential equations (it is a
tau-function of the Korteweg--de Vries integrable hierarchy obeying
the ``string equation'') which are sufficient to compute all the
intersection numbers explicitly. One way to write the equations is
as Virasoro conditions
\[
L_k(e^F)=0,\qquad k=-1,0,1,2\dots,
\]
for certain differential operators $L_k$ of order at most 2 obeying
the commutation relations $[L_{j},L_{k}]=(j-k)L_{j+k}$ of the Lie
algebra of polynomial vector fields.

Before Okounkov  few results were available for general projective
varieties $V$ and they were mostly restricted to genus $g=0$
Gromov--Witten invariants (quantum cohomology). For our purpose the
conjecture of Eguchi, Hori and Xiong is relevant here. Again,
Gromov--Witten invariants of $V$ can be encoded into a generating
function $F_V$ depending on variables $t_{j,a}$ where $a$ labels a
basis of the cohomology of $V$. Eguchi, Hori and Xiong extended
Witten's definition of the differential operators $L_k$ and
conjectured that $F_V$ obeys the Virasoro conditions
$L_k(e^{F_V})=0$ with these operators.

\section{Gromov--Witten invariants of curves}

In a remarkable series of papers \cites{OP1,OP2,OP3}, Okounkov and
Pandharipande give an exhaustive description of the Gromov--Witten
invariants of curves. They prove the Eguchi--Hori--Xiong conjecture
for general projective curves $V$, give explicit descriptions in the
case of genus 0 and 1, show that the generating function for
$V=\mathbb{P}^1$ is a tau-function of the Toda hierarchy and
consider also in this case the $\mathbb C^\times$-equivariant
theory, which is shown to be governed by the 2D-Toda hierarchy. They
also show that GW invariants of $V=\mathbb P^1$ are unexpectedly
simple and more basic than the GW invariants of a point, in the
sense that the latter can be obtained as a limit, giving thus a more
transparent proof of Kontsevich's theorem.

 A key ingredient is the {\em Gromov--Witten/Hurwitz
correspondence} relating GW invariants of a curve $V$ to {\em
Hurwitz numbers}, the numbers of branched covering of $V$ with given
ramification type at given points. A basic beautiful formula of
Okounkov and Pandharipande is the formula for the {\em stationary}
GW invariants of a curve $V$ of genus $g(V)$, namely those for the
Poincar\'e dual $\omega$ of a point:
\begin{equation}\label{e-GW-for-curves}
\langle\tau_{k_1}(\omega) \cdots
\tau_{k_n}(\omega)\rangle_{\beta=d\cdot [V],g}^{\bullet\, V} =
\sum_{|\lambda|=d}
\left(\frac{\mathrm{dim}\,\lambda}{d!}\right)^{2-2g(V)}
\prod_{i=1}^{n} \frac{p_{k_i+1}(\lambda)}{(k_i+1)!}.
\end{equation}
The (finite!) summation is over all partitions $\lambda$ of the
degree $d$ and $\mathrm{dim}\,\lambda$ is the dimension of the
corresponding irreducible representation of $S_d$. The genus $g$ of
the domain is fixed by the condition that the cohomological degree
of the integrand is equal to the dimension of the virtual
fundamental class. It is convenient here to include also stable maps
with possibly disconnected domains and this is indicated by the
bullet. The functions $p_k(\lambda)$ on partitions are described
below.

Hurwitz numbers can be computed combinatorially and are given in
terms of representation theory of the symmetric group by an explicit
formula of Burnside. If the covering map at the $i$th point looks
like $z\to z^{k_i+1}$, i.e., if the monodromy at the $i$th point is
a cycle of length $k_i+1$, the formula is
\[
H^V_d(k_1+1,\dots,k_n+1) = \sum_{|\lambda|=d}
\left(\frac{\mathrm{dim}\,\lambda}{d!}\right)^{2-2g(V)}
\prod_{i=1}^{n} {f_{k_i+1}(\lambda)}.
\]
Thus in this case the GW/Hurwitz correspondence is given by the
substitution rule $f_{k+1}(\lambda)\to p_{k+1}(\lambda)/(k+1)!$. The
functions $f_k$ and $p_k$ are basic examples of {\em shifted
symmetric functions}, a theory initiated by Kerov and Olshanski, and
the results of Okounkov and Pandharipande offer a geometric
realization of this theory. A shifted symmetric polynomial of $n$
variables $\lambda_1,\dots,\lambda_n$ is a polynomial invariant
under the action of the symmetric group given by permuting
$\lambda_j-j$. A shifted symmetric function is a function of
infinitely many variables $\lambda_1,\lambda_2\dots,$ restricting
for each $n$ to a shifted symmetric polynomial of $n$ variables if
all but the first $n$ variables are set to zero. Shifted symmetric
functions form an algebra $\Lambda^*=\mathbb Q[p_1,p_2,\dots]$
freely generated by the {\em regularized shifted power sums},
appearing in the GW invariants:
\[
p_k(\lambda)=\sum_j\left(
\left(\lambda_j-j+\textstyle{\frac12}\right)^k-\left(-j+\textstyle{\frac12}\right)^k\right)
+(1-2^{-k})\zeta(-k)
\]
The second  term and the Riemann zeta value ``cancel out'' in the
spirit of Ramanujan's second letter to Hardy:
$1+2+3+\cdots=-\frac1{12}$. The shifted symmetric functions
$f_k(\lambda)$ appearing in the Hurwitz numbers are central
characters of the symmetric groups $S_n$:
$f_1=|\lambda|=\sum\lambda_i$ and the sum of the elements of the
conjugacy class of a cycle of length $k \geq 2$ in the symmetric
group $S_n$ is a central element acting as $f_k(\lambda)$ times the
identity in the irreducible representation corresponding to
$\lambda$. The functions $p_k$ and $f_k$ are two natural shifted
versions of Newton power sums.
 \dontprint{ More generally, one has the map sending a conjugacy
class $\eta$, given by its cycle lengths
$\eta_1\geq\cdots\geq\eta_k>1$, to its eigenvalues $f_\eta(\lambda)$
in irreducible representations. These functions form a basis of
$\Lambda^*$ by a theorem of Kerov--Olshanski. The expression of
$p_k/k!$ in terms of the basis $f_\eta$ starts with $f_k$ and has
positive coefficients (except for the constant term with
$\eta=\varnothing$). Thus one can formulate the GW/Hurwitz
correspondence as the statement that GW-invariants are obtained from
Hurwitz number by replacing $k$-cycles by {\em completed $k$-cycles}
(a notion due to Okounkov--Pandharipande--Zorich), which are
$k$-cycles corrected by a universal positive linear combination of
(smaller) permutations and a constant term.
 }

In the case of genus $g(V)=0,1$ Okounkov and Pandharipande
reformulate \eqref{e-GW-for-curves} in terms of expectation values
and traces in fermionic Fock spaces and get more explicit
descriptions and recursion relations. In particular if $V=E$ is an
elliptic curve the generating function of GW invariants reduces to
the formula of Bloch and Okounkov \cite{BO} for the character of the
infinite wedge projective representation of the algebra of
polynomial differential operators, which is expressed in terms of
Jacobi theta functions. As a corollary, one obtains that
\[
\sum_{d}q^d
\langle\tau_{k_1}(\omega)
\cdots
\tau_{k_n}(\omega)\rangle_{d}^{\bullet\, E}
\]
belongs to the ring $\mathbb Q[E_2,E_4,E_6]$ of {\em quasimodular
forms}. \dontprint{If all $k_i=1$ these GW invariants are Hurwitz
numbers for simple ramifications (since $f_2=p_2$) and one recovers
a result of Dijkgraaf and Kaneko--Zagier.}

A shown by Eskin and Okounkov \cite{EO}, one can use quasimodularity
to compute the asymptotics as $d\to \infty$ of the number of
connected ramified degree $d$ coverings of a torus with given
monodromy at the ramification points. By a theorem of
Kontsevich--Zorich and Eskin--Mazur, this asymptotics gives the
volume of the moduli space of holomorphic differentials on a curve
with given orders of zeros, which is in turn related to the dynamics
of billiards in rational polygons. Eskin and Okounkov give explicit
formulae for these volumes and prove in particular the
Kontsevich--Zorich conjecture that they belong to $\pi^{-2g}\mathbb
Q$ for curves of genus $g$.

\section{Donaldson--Thomas invariants}\label{s-DT}
As is clear from the dimension formula \eqref{e-1} the case of
three-dimensional varieties $V$ plays a very special role. In this
case, which in the Calabi-Yau case $K_V=0$ is the original context
studied in string theory, it is possible to define invariants
counting curves by describing curves by equations rather than in
parametric form. Curves in $V$ of genus $g$ and class $\beta\in
H_2(V)$ given by equations are parametrized by Grothendieck's
Hilbert schemes $\mathrm{Hilb}(V;\beta,\chi)$ of subschemes of $V$
with given Hilbert polynomial of degree 1. The invariants $\beta$,
$\chi=2-g$ are encoded in the coefficients of the Hilbert
polynomial. R. Thomas constructed a virtual fundamental class of
$\mathrm{Hilb}(V; \beta,\chi)$ for three-folds $V$ of dimension
$-\beta\cdot K_V$, the same as the dimension of $\overline
M_{g,0}(V,\beta)$. Thus one can define Donaldson--Thomas (DT)
invariants as intersection numbers on this Hilbert scheme. There is
no direct geometric relation between $\mathrm{Hilb}(V;\beta,\chi)$
and $\overline M_{g,0}(V,\beta)$, and indeed the (conjectural)
relation between Gromov--Witten invariants and Donaldson--Thomas
invariants is quite subtle. In its simplest form, it relates the GW
invariants $ \int_{\bar
M^\bullet_{g,n}(V,\beta)}\prod\mathrm{ev}_i^*\gamma_i$ to the DT
invariants $ \int_{\mathrm{Hilb}(V;\beta,\chi)}\prod c_2(\gamma_i)$.
The class $c_2(\gamma)$ is the coefficient of $\gamma\in H^*(V)$ in
the K\"unneth decomposition of the second Chern class of the ideal
sheaf of the universal family $\mathcal V\subset
\mathrm{Hilb}(V;\beta,\chi)\times V$.

The conjecture of Maulik, Nekrasov, Okounkov and Pandharipande
\cites{MNOP1, MNOP2}, inspired by ideas of string theory \cite{ORV}
states that suitably normalized generating functions
$Z'_{GW}(\gamma;u)_\beta$, $Z'_{DT}(\gamma;q)_\beta$ are essentially
related by a coordinate transformation:
\[
(-iu)^{-d}Z'_{GW}(\gamma_1,\dots,\gamma_n;u)_\beta
=q^{-d/2}Z'_{DT}(\gamma_1,\dots,\gamma_n;q)_\beta,
\quad \text{if $q=-e^{iu}$, $\beta\neq 0$}.
 \]
Here $d=-\beta\cdot K_V$ is the virtual dimension. Moreover these
authors conjecture that there $Z'_{DT}(\gamma;q)_\beta$ is a {\em
rational} function of $q$. This has the important consequence that
all (infinitely many) GW invariants are determined in principle by
finitely many DT invariants. Versions of these conjectures are
proven for local curves and the total space of the canonical bundle
of a toric surface. The GW/DT correspondence can be viewed as a
far-reaching generalization of formula \eqref{e-GW-for-curves}, to
which it reduces in the case where $V$ is the product of a curve
with $\mathbb C^2$.

\section{Other uses of partitions}
Here is a short account of other results of Okounkov based on the
occurrence of partitions.

One early result of Okounkov \cite{O} is his first proof of the
Baik--Deift--Johansson conjecture (two further different proofs
followed, one by Borodin, Okounkov and Olshanski and one by
Johansson). This conjecture states that, as $n\to\infty$, the joint
distribution of the first few rows of a random partition of $n$ with
the Plancherel measure
$P(\lambda)=(\mathrm{dim}\,\lambda)^2/|\lambda|!$, natural from
representation theory, is the same, after proper shift and
rescaling, as the distribution of the first few eigenvalues of a
Gaussian random hermitian matrix of size $n$. The proof involves
comparing random surfaces given by Feynman diagrams and by ramified
coverings and contains many ideas that anticipate Okounkov's later
work on Gromov--Witten invariants.

Random partitions also play a key role in the work \cite{NO} of
Nekrasov and Okounkov on $N=2$ supersymmetric gauge theory in four
dimensions. Seiberg and Witten gave a formula for the effective
``prepotential'', postulating a duality with a theory of monopoles.
The Seiberg--Witten formula is given in terms of periods on a family
of algebraic curves, closely connected with classical integrable
systems. Nekrasov showed how to rigorously define the prepotential
of the gauge theory as a regularized instanton sum given by a
localization integral on the moduli space of antiselfdual
connections on $\mathbb R^4$. Nekrasov and Okounkov show that this
localization integral can be written in terms of a measure on
partitions with periodic potential and identify the Seiberg--Witten
prepotential with the surface tension of the limit shape.

Partitions of $n$ also label $(\mathbb C^\times)^2$-invariant ideals
of codimension $n$ in $\mathbb C[x,y]$ and thus appear in
localization integrals on the Hilbert scheme of points in the plane.
Okounkov and Pandharipande \cite{OP} describe the ring structure of
the equivariant quantum cohomology (genus zero GW invariants) of
this Hilbert scheme in terms of a time-dependent version of the
Calogero--Moser operator from integrable systems.

\section{Dimers}
Dimers are a much studied classical subject in statistical mechanics
and graph combinatorics. Recent spectacular progress in this subject
is due to the discovery by Okounkov and collaborators of a close
connection of planar dimer models with real algebraic geometry.
\dontprint{In particular, a complete description of the phase
diagram and of the asymptotic shape of interfaces is obtained.}

A {\em dimer configuration} (or perfect matching) on a bipartite
graph $G$ is a subset of the set of edges of $G$ meeting every
vertex exactly once. For example if $G$ is a square grid we may
visualize a dimer configuration as a tiling of a checkerboard by
dominoes. In statistical mechanics one assigns positive weights
(Boltzmann weights) to edges of $G$ and defines the weight of a
dimer configuration as the product of the weights of its edges. The
basic tool is the Kasteleyn matrix of $G$, which is up to certain
signs the weighted adjacency matrix of $G$. For finite $G$ Kasteleyn
proved that the partition function (i.e., the sum of the weights of
all dimer configurations) is the absolute value of the determinant
of the Kasteleyn matrix.

Kenyon, Okounkov and Sheffield consider a doubly periodic bipartite
graph $G$ embedded in the plane with doubly periodic weights. For
each natural number $n$ one then has a probability measure on dimer
configurations on $G_n= G/n\mathbb Z^2$ and statistical mechanics of
dimers is essentially the study of the asymptotics of these
probability measures in the thermodynamic limit $n\to\infty$. One
key observation is the Kasteleyn matrix on $G_1$ can be twisted by a
character $(z,w)\in(\mathbb C^\times)^2$ of $\mathbb Z^2$ and thus
one defines the {\em spectral curve} as the zero set $P(z,w)=0$ of
the determinant of the twisted Kasteleyn matrix $P(z,w)=\det
K(z,w)$. This determinant is a polynomial in $z^{\pm1}, w^{\pm1}$
with real coefficients and thus defines a real plane curve.

The main observation of Kenyon, Okounkov and Sheffield \cite{KOS} is
that the spectral curve belongs to the very special class of
(simple) Harnack curves, which were studied in the 19th century and
have reappeared recently in real algebraic geometry. Kenyon,
Okounkov and Sheffield show that in the thermodynamic limit, three
different {\em phases} (called gaseous, liquid and frozen) arise.
These phases are characterized by qualitatively different
long-distance behaviour of pair correlation functions. One can see
these phases by varying two real parameters $(B_1,B_2)$ ( the
``magnetic field'') in the weights, so that the spectral curve
varies by rescaling the variables. The regions in the
$(B_1,B_2)$-plane corresponding to different phases are described in
terms of the {\em amoeba} of the spectral curve, namely the image of
the curve by the map
$\mathrm{Log}\colon(z,w)\mapsto(\log|z|,\log|w|)$. The amoeba of a
curve is a closed subset of the plane which looks a bit like the
microorganism with the same name. The amoeba itself corresponds to
the liquid phase, the bounded components of its complement to the
gaseous phase and the unbounded components to the frozen phase. This
insight has a lot of consequences for the statistics of dimer models
and lead Okounkov and collaborators to beautiful results on
interfaces with various boundary conditions \cites{KOS,KO3}.

Such a precise and complete description of phase diagrams and shapes
of interfaces is unprecedented in statistical mechanics.

\section{Random surfaces}
One useful interpretation of dimers is as models for random surfaces
in three-dimensional space. In the simplest case one considers a
model for a melting or dissolving cubic crystal in which at a corner
some atoms are missing (Fig.~1).

\centerline{
\begin{picture}(150,130)(0,0)
\scalebox{0.3}{\includegraphics{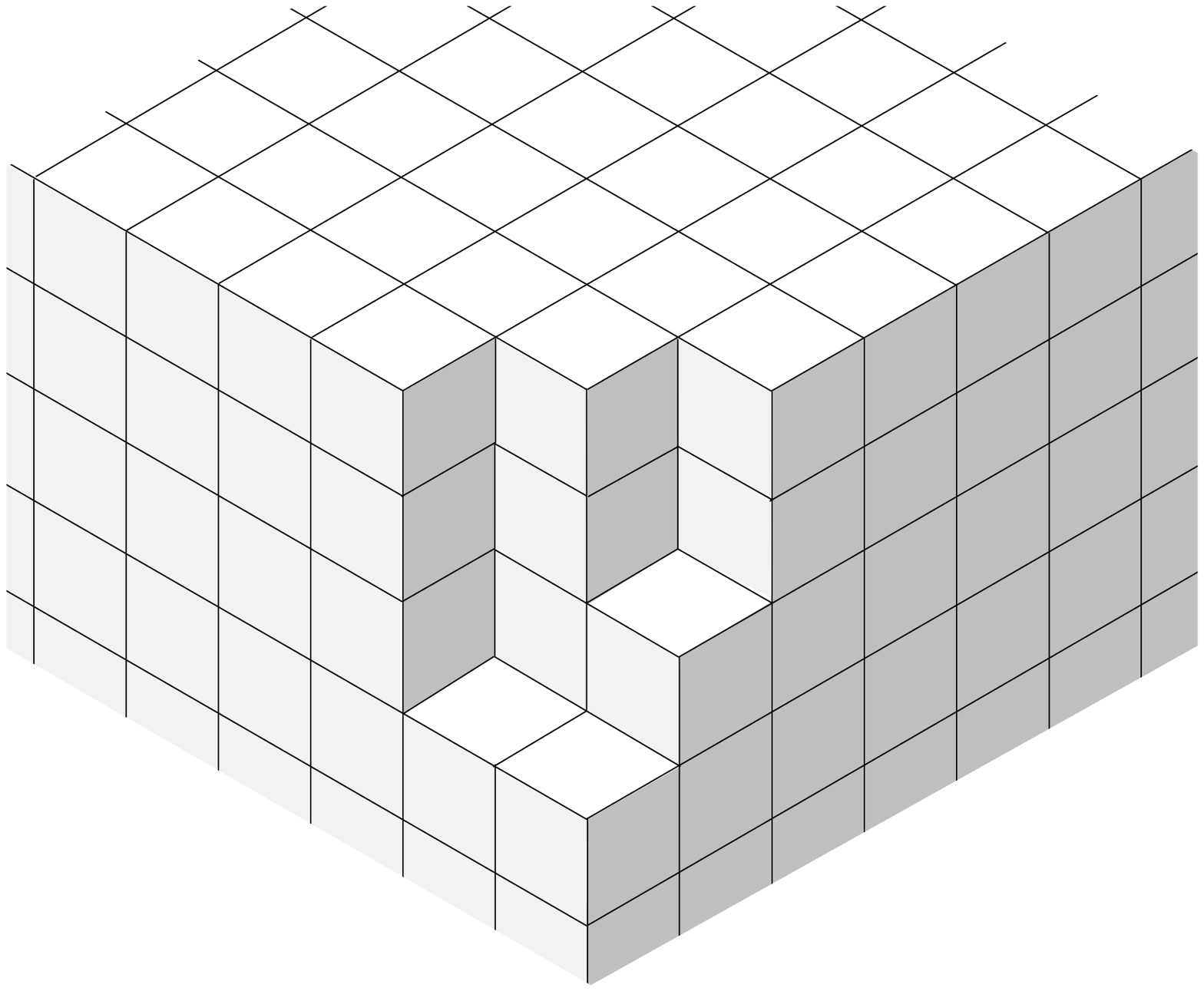}}
\end{picture} \ \ \ \ \ \ \ \ \ \ \ \
\begin{picture}(130,130)(0,-20)
\scalebox{0.3}{\includegraphics{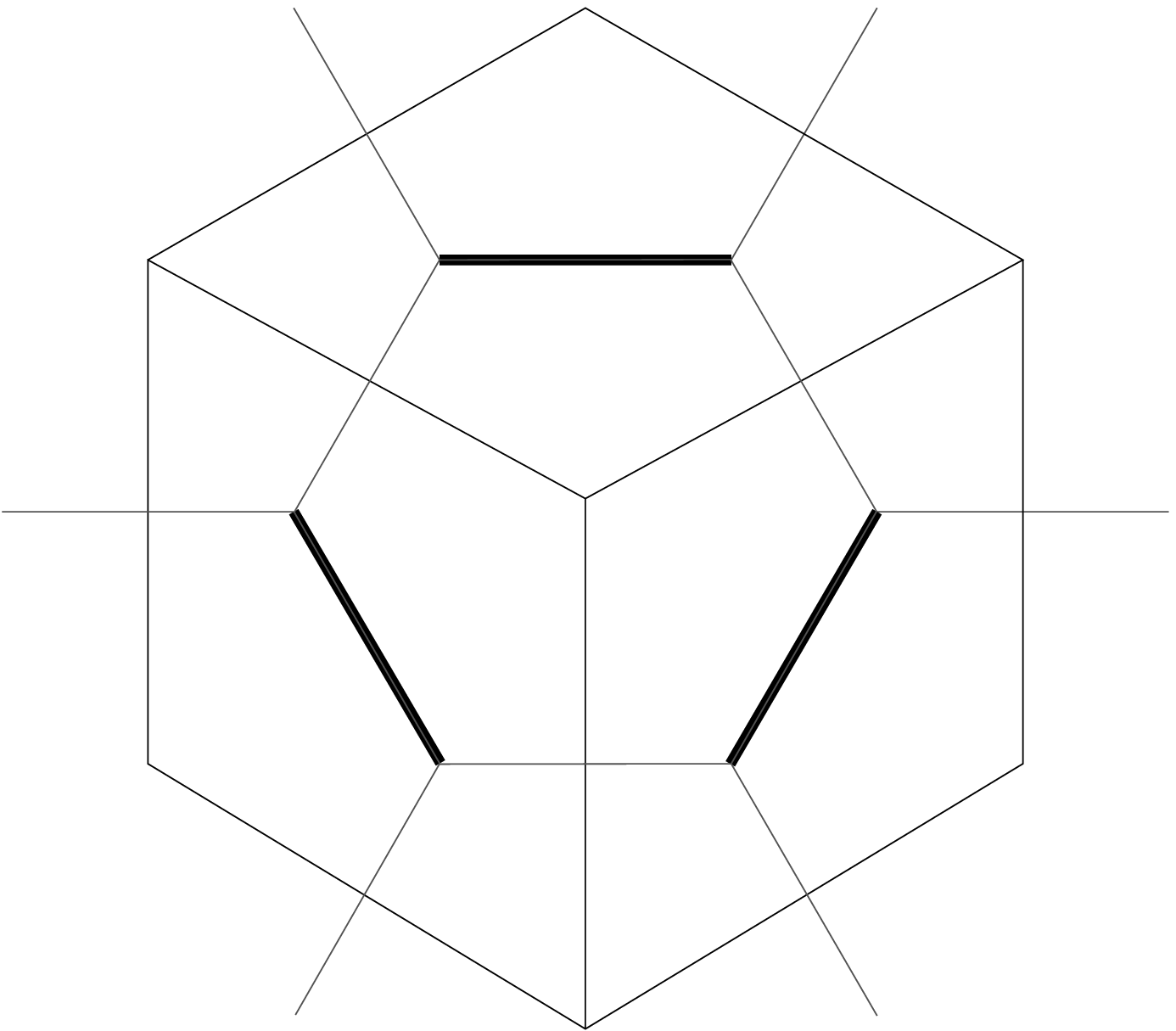}}
\end{picture}
} \centerline{\small {\sc Fig.~1.} A melting crystal corner (left)
and the relation between tilings and dimers (right)}
\medskip

\noindent Viewing the corner from the (1,1,1) direction one sees a
tiling of the plane by $60^\circ$ rhombi, which is the same as a
dimer configuration on a honeycomb lattice (each tile covers one
dimer of the dimer configuration). In this simple model one gives
the same probability for every configuration with given missing
volume. If one lets the size of the cubes go to zero keeping the
missing volume fixed, the probability measure concentrates on an a
surface, the limit shape. More generally, every planar dimer model
can be rephrased as a random surface model and limit shapes for more
general crystal corner geometries can be defined. Kenyon, Okounkov
and Sheffield show that the limit shape is given by the graph of
(minus) the {\em Ronkin function} $R(x,y)=(2\pi i)^{-2}
\int_{|z|=|w|=1} \log(P(e^xz,e^yw))dzdw/zw$ of the spectral curve
(in the case of the honeycomb lattice with equal weights,
$P(z,w)=z+w+1$). This function is affine on the complement of the
amoeba and strictly convex on the amoeba. So the connected
components of the complement of the amoeba are the projections of
the facets of the melting crystal.

In addition to this surprising connection with real algebraic
geometry, random surfaces of this type are essential in the GW/DT
correspondence, see Section \ref{s-DT}, as they arise in
localization integrals for DT invariants of toric varieties.

\section{The moduli space of Harnack curves}
The notions used by Okounkov and collaborators in their study of
dimer models arose in an independent recent development in real
algebraic geometry. Their result bring a new probabilistic point of
view in this classical subject.

In real algebraic geometry, unsurmountable difficulties already
appear when one consider curves. The basic open question is the
first part of Hilbert's 16th problem: what are the possible
topological types of a smooth curve in the plane given by a
polynomial equation $P(z,w)=0$ of degree $d$?\dontprint{ Harnack
proved that such a curve has at most $1+(d-1)(d-2)/2$  connected
components (``ovals''). Curves with this many ovals are called
maximal.} Topological types up to degree 7 are known but very few
general results are available. In a recent development in real
algebraic geometry in the context of toric varieties the class of
Harnack curves plays an important role and can be characterized in
many equivalent way. In one definition, due to Mikhalkin, a Harnack
curve is a curve such that the map to its amoeba is 2:1 over the
interior, except at possible nodal points; equivalently, by a
theorem of Mikhalkin and Rullg{\aa}rd, a Harnack curve is a curve
whose amoeba has area equal to the area of the Newton polygon of the
polynomial $P$. These equivalent properties determine the
topological type completely.

Kenyon and Okounkov prove \cite{KO2} that {\em every} Harnack curve
is the spectral curve of some dimer model. They obtain an explicit
parametrization of the moduli space of Harnack curves with fixed
Newton polygon by weights of dimer models, and deduce in particular
that the moduli spaces are connected.

\section{Concluding remarks}
Andrei Okounkov is a highly creative mathematician with both an
exceptional breadth and a sense of unity of mathematics, allowing
him to use and develop, with perfect ease, techniques and ideas from
all branches of mathematics to reach his research objectives. His
results not only settle important questions and open new avenues of
research in several fields of mathematics, but they have the
distinctive feature of mathematics of the very best quality: they
give simple complete answers to important natural questions, they
reveal hidden structures and new connections between mathematical
objects and they involve new ideas and techniques with wide
applicability.

\dontprint{
 Namely:
1.~They give simple, unexpected and complete answers to natural
questions---for example, given the complicated structure of GW
invariants for a point, it came as a surprise that GW invariants of
curves are given by the {\em finite} sum of
eq.~\eqref{e-GW-for-curves} and that they could be computed in full
generality. 2.~They reveal hidden structures and new connections
between mathematical objects---such as the relation of GW invariants
with shifted symmetric function or the role of Harnack curves in the
statistical mechanics of dimers. 3.~They involve new ideas and
techniques that can be applied to other problems---for example,
Okounkov showed us that asymptotics of random partitions are
significant in the most diverse fields of mathematics. }

Moreover, in addition to obtaining several results of this quality
representing significant progress in different fields, Okounokov is
able to create the ground, made of visions, intuitive ideas and
techniques, where new mathematics appears. A striking example for
this concerns the relation to physics:  many important developments
in mathematics of the last few decades have been inspired by high
energy physics, whose intuition is based on notions often
inaccessible to mathematics. Okounkov's way of proceeding is to
develop a mathematical intuition alternative to the intuition of
high energy physics, allowing him and his collaborators to go beyond
the mere verification of predictions of physicists. Thus, for
example, in approaching the topological vertex of string theory,
instead of stacks of D-branes and low energy effective actions we
find mathematically more familiar notions such as localization and
asymptotics of probability measures. As a consequence, the scope of
Okounkov's research programme goes beyond the context suggested by
physics: for example the Maulik--Nekrasov--Okounkov--Pandharipande
conjecture is formulated (and proved in many cases) in a setting
which is much more general than the Calabi--Yau case arising in
string theory.


\begin{bibdiv}
\begin{biblist}

\bib{BO}{article}{
   author={Bloch, Spencer},
   author={Okounkov, Andrei},
   title={The character of the infinite wedge representation},
   journal={Adv. Math.},
   volume={149},
   date={2000},
   number={1},
   pages={1--60},
   issn={0001-8708},
   review={\MR{1742353 (2001g:11059)}},
}

\bib{EO}{article}{
   author={Eskin, Alex},
   author={Okounkov, Andrei},
   title={Asymptotics of numbers of branched coverings of a torus and
   volumes of moduli spaces of holomorphic differentials},
   journal={Invent. Math.},
   volume={145},
   date={2001},
   number={1},
   pages={59--103},
   issn={0020-9910},
   review={\MR{1839286 (2002g:32018)}},
}

\bib{KOS}{article}{
   author={Kenyon, Richard},
   author={Okounkov, Andrei},
   author={Sheffield, Scott},
   title={Dimers and amoebae},
   journal={Ann. of Math. (2)},
   volume={163},
   date={2006},
   number={3},
   pages={1019--1056},
   issn={0003-486X},
   review={\MR{2215138}},
}
\bib{KO2}{article}{
   author={Kenyon, Richard},
   author={Okounkov, Andrei},
   title={Planar dimers and Harnack curves},
   journal={Duke Math. J.},
   volume={131},
   date={2006},
   number={3},
   pages={499--524},
   issn={0012-7094},
   review={\MR{2219249}},
}
\bib{KO3}{article}{
    title = {{Limit shapes and the complex burgers equation}},
    author = {Richard Kenyon and Andrei Okounkov},
    eprint = {arXiv:math-ph/0507007}
}
\bib{MNOP1}{article}{
    title = {{Gromov-Witten theory and Donaldson-Thomas theory, I}},
    author = {Maulik, D.},
    author = {Nekrasov, N.},
    author = {Okounkov, A.},
    author = {Pandharipande, R.},
    eprint = {arXiv:math.AG/0312059}}
\bib{MNOP2}{article}{
    title = {{Gromov-Witten theory and Donaldson-Thomas theory, II}},
    author = {Maulik, D.},
    author = {Nekrasov, N.},
    author = {Okounkov, A.},
    author = {Pandharipande, R.},
    eprint = {arXiv:math.AG/0406092}}

\bib{NO}{article}{
   author={Nekrasov, Nikita A.},
   author={Okounkov, Andrei},
   title={Seiberg-Witten theory and random partitions},
   conference={
      title={The unity of mathematics},
   },
   book={
      series={Progr. Math.},
      volume={244},
      publisher={Birkh\"auser Boston},
      place={Boston, MA},
   },
   date={2006},
   pages={525--596},
   review={\MR{2181816}},
}
\bib{O}{article}{
   author={Okounkov, Andrei},
   title={Random matrices and random permutations},
   journal={Internat. Math. Res. Notices},
   date={2000},
   number={20},
   pages={1043--1095},
   issn={1073-7928},
   review={\MR{1802530 (2002c:15045)}},
}
\bib{OP1}{article}{
   author={Okounkov, Andrei},
   author={Pandharipande, Rahul},
   title={Gromov-Witten theory, Hurwitz theory, and completed cycles},
   journal={Ann. of Math. (2)},
   volume={163},
   date={2006},
   number={2},
   pages={517--560},
   issn={0003-486X},
   review={\MR{2199225}},
}
\bib{OP2}{article}{
   author={Okounkov, Andrei},
   author={Pandharipande, Rahul},
   title={The equivariant Gromov-Witten theory of ${\bf P}\sp 1$},
   journal={Ann. of Math. (2)},
   volume={163},
   date={2006},
   number={2},
   pages={561--605},
   issn={0003-486X},
   review={\MR{2199226 (2006j:14075)}},
}
\bib{OP3}{article}{
   author={Okounkov, Andrei},
   author={Pandharipande, Rahul},
   title={Virasoro constraints for target curves},
   journal={Invent. Math.},
   volume={163},
   date={2006},
   number={1},
   pages={47--108},
   issn={0020-9910},
   review={\MR{2208418}},
}

\bib{OP}{article}{
   title = {{Quantum cohomology of the Hilbert scheme of points in the
        plane}},
    author = {Andrei Okounkov and Rahul Pandharipande},
    eprint = {arXiv:math.AG/0411210}}

\bib{ORV}{article}{
   author={Okounkov, Andrei},
   author={Reshetikhin, Nikolai},
   author={Vafa, Cumrun},
   title={Quantum Calabi-Yau and classical crystals},
   conference={
      title={The unity of mathematics},
   },
   book={
      series={Progr. Math.},
      volume={244},
      publisher={Birkh\"auser Boston},
      place={Boston, MA},
   },
   date={2006},
   pages={597--618},
   review={\MR{2181817}},
}
\end{biblist}
\end{bibdiv}
\end{document}